\documentclass{article}
\usepackage[utf8]{inputenc}
\usepackage{amsmath}
\usepackage[T1]{fontenc}    
\usepackage{url}            
\usepackage{booktabs}       
\usepackage{amsfonts}       
\usepackage{nicefrac}       
\usepackage{microtype}      
\usepackage{lipsum}
\usepackage{graphicx}
\usepackage[english]{babel}
\usepackage{amsthm}
\usepackage{tabularx}
\usepackage{array}
\usepackage{mathrsfs}
\usepackage{mathtools}
\usepackage{amssymb}

\newtheorem{theorem}{Theorem}[section]
\newtheorem{conjecture}{Conjecture}[section]
\newtheorem{remark}{Remark}[section]
\newtheorem{corollary}{Corollary}[theorem]
\newtheorem{lemma}[theorem]{Lemma}

\begin{document}

\title{On the quantity $m^2 - p^k$ where $p^k m^2$ \\ is an odd perfect number - Part II \footnote{2010 \it{Mathematics Subject Classification}. 11A05, 11A25. \\
\it{Key words and phrases}. Odd perfect number, Descartes-Frenicle-Sorli Conjecture, Dris Conjecture. \\
\it{Corresponding author}. Jose Arnaldo Bebita Dris}}
\author{Jose Arnaldo Bebita Dris$^1$ and Immanuel Tobias San Diego$^2$ \\ 
  \texttt{$^1$ Graduate of De La Salle University, Manila, Philippines} \\
	\texttt{$^1$ josearnaldobdris@gmail.com} \\
	\texttt{$^2$ Trinity University of Asia, Quezon City, Philippines} \\
	\texttt{$^2$ itsandiego@tua.edu.ph} \\}

\maketitle

\begin{abstract}
Let $p^k m^2$ be an odd perfect number with special prime $p$.  Extending previous work of the authors, we prove that the inequality $m < p^k$ follows from $m^2 - p^k = 2^r t$, where $r \geq 2$ and $\gcd(2,t)=1$, under the following hypotheses:
\begin{enumerate}
\item{$m > t > 2^r$; or}
\item{$m > 2^r > t$.}
\end{enumerate}
We also prove that the estimate $m^2 - p^k > 2m$ holds.  We can also improve this unconditional estimate to $m^2 - p^k > {313m^2}/315$.
\end{abstract}

\section{Introduction}\label{Intro}

Let $\sigma(x)$ be the sum of the divisors of the positive integer $x$.  Denote the deficiency [\url{https://oeis.org/A033879}] of $x$ by $D(x)=2x-\sigma(x)$, and the aliquot sum [\url{https://oeis.org/A001065}] of $x$ by $s(x)=\sigma(x)-x$.  Note that we have the identity $$D(x) + s(x) = x.$$

If a positive integer $N$ is odd and $\sigma(N)=2N$, then $N$ is called an odd perfect number.  Euler proved that a hypothetical odd perfect number, if one exists, must have the so-called \textit{Eulerian form} $$N = p^k m^2,$$ where $p$ is the \textit{special prime} satisfying $p \equiv k \equiv 1 \pmod 4$ and $\gcd(p,m)=1$.  Despite extensive computer searches, to date nobody has found an odd perfect number.  Ochem and Rao \cite{OchemRao} has proved that $N > {10}^{1500}$, making the existence of odd perfect numbers appear very unlikely.

Descartes, Frenicle, and subsequently Sorli conjectured that $k=1$ always holds \cite{Beasley}.  Sorli predicted that $k=1$ is true after testing large numbers with eight distinct prime factors for perfection \cite{Sorli}.  Dris \cite{Dris2}, and Dris and Tejada (\cite{DrisTejada},\cite{DrisTejada2}), call this conjecture as the Descartes--Frenicle--Sorli Conjecture, and derive conditions equivalent to $k=1$.

Dris conjectured in \cite{Dris} that the factors $p^k$ and $m$ are related by the inequality $p^k < m$.  Brown was the first to show that the inequality $p < m$ holds in a preprint \cite{Brown}.  Other proofs for the estimate $p < m$ have been given by Dris \cite{Dris3}, Starni \cite{Starni}, and Dagal and Dris \cite{DagalDris2}.  (Note that if one could disprove the Dris Conjecture, so that one would have $p < m < p^k$, it would follow that the Descartes--Frenicle--Sorli Conjecture is false.)

Since $m$ is odd, then $m^2 \equiv 1 \pmod 4$.  Likewise, $p \equiv k \equiv 1 \pmod 4$ is true, which implies that $p^k \equiv 1 \pmod 4$ holds.  It follows that $m^2 - p^k \equiv 0 \pmod 4$.  Since
$$p^k < \frac{2m^2}{3}$$
(by a result of Dris~\cite{Dris}), we know \textit{a priori} that
$$m^2 - p^k > \frac{p^k}{2}$$
so that we are sure that $m^2 - p^k > 0$.  In particular, since $m^2 - p^k \equiv 0 \pmod 4$, we infer that $m^2 - p^k \geq 4$.

The index $i(p)$ of the odd perfect number $N = p^k m^2$ at the prime $p$ is then equal to
$$i(p):=\frac{\sigma(N/p^k)}{p^k}=\frac{\sigma(m^2)}{p^k}=\frac{m^2}{\sigma(p^k)/2}=\frac{D(m^2)}{s(p^k)}=\frac{s(m^2)}{D(p^k)/2}=\gcd(m^2,\sigma(m^2)).$$
The term \emph{index of an odd perfect number (at a certain prime)} was coined by Chen and Chen \cite{ChenChen}.

The following lemma gives a nontrivial lower bound for $i(p)$, which was proved by Broughan, Delbourgo, and Zhou \cite{BroughanDelbourgoZhou}.

\begin{lemma}\label{BroughanDelbourgoZhou}
If $N = p^k m^2$ is an odd perfect number given in Eulerian form, then $i(p) \geq 315$.
\end{lemma}

Proving the following lemma is trivial, and this follows from the estimate $p^k < m^2$ as proved in Dris \cite{Dris}, and the lower bound for the magnitude of an odd perfect number as proved in Ochem and Rao \cite{OchemRao}:

\begin{lemma}\label{estimateform}
If $N = p^k m^2$ is an odd perfect number given in Eulerian form, then $$m > \sqrt[4]{N} > {10}^{375}.$$
\end{lemma}

Finally, recall that we obtained the following results from an earlier paper of the authors \cite{DrisSanDiego} on this topic:

\begin{theorem}\label{QuantityIsNotASquare}
If $N = p^k m^2$ is an odd perfect number given in Eulerian form, then $m^2 - p^k$ is not a square.
\end{theorem}

\begin{theorem}\label{PreludeToMainResult}
If $N = p^k m^2$ is an odd perfect number given in Eulerian form and satisfying $m^2 - p^k = 8$, then the inequality $m < p^k$ holds.
\end{theorem}

In the present paper, we present some general conditions under which the estimate $m < p^k$ is true.

\section{The results}\label{Summary}

We now present a summary of our results in this section.

The first proposition gives us a very large lower bound for the quantity $m^2 - p^k$.

\begin{theorem}\label{VeryLargeLowerBound}
If $N = p^k m^2$ is an odd perfect number given in Eulerian form, then $m^2 - p^k > 2m$.
\end{theorem}

We can also prove the following corollary, which modestly improves on Theorem \ref{VeryLargeLowerBound}.

\begin{corollary}\label{VeryLargeLowerBound2}
If $N = p^k m^2$ is an odd perfect number given in Eulerian form, then $$m^2 - p^k > \frac{313m^2}{315}.$$
\end{corollary}

Next, in the second and third propositions, we derive some general conditions under which the inequality $m < p^k$ holds.

\begin{theorem}\label{MainResult1}
Let $N = p^k m^2$ be an odd perfect number given in Eulerian form and satisfying $m^2 - p^k = 2^r t$, where $r \geq 2$ and $\gcd(2,t)=1$.  If $2^r t < (m - 1)^2$, then $m < p^k$.
\end{theorem}

\begin{theorem}\label{MainResult2}
Let $N = p^k m^2$ be an odd perfect number given in Eulerian form and satisfying $m^2 - p^k = 2^r t$, where $r \geq 2$ and $\gcd(2,t)=1$.
\begin{enumerate}
\item{If $m > t > 2^r$, then the inequality $m < p^k$ holds.}
\item{If $m > 2^r > t$, then the inequality $m < p^k$ holds.}
\end{enumerate}
\end{theorem}

\section{A proof of Theorem \ref{VeryLargeLowerBound}}\label{Proof1}

Let $N = p^k m^2$ be an odd perfect number with special prime $p$.  Assume to the contrary that $m^2 - p^k \leq 2m$.

Since $\gcd(p,m)=1$, then we can consider two cases:

(1) Suppose that $p^k < m$.  By assumption, we have $m^2 \leq p^k + 2m$.  This implies that $$m^2 < m + 2m = 3m.$$
This gives $m < 3$, which contradicts Lemma \ref{estimateform}.

(2) Suppose that $m < p^k$.  By assumption, we have $m^2 \leq p^k + 2m$.  But the inequality $m < p^k$ together with the inequality $m^2 \leq p^k + 2m$ will contradict the lower bound $\sigma(m^2)/p^k \geq 7$ by Dris and Luca \cite{DrisLuca}, as follows:
$$\Bigg((m^2 \leq p^k + 2m) \land (m < p^k)\Bigg) \implies (m^2 < 3p^k).$$
However, the estimate
$$m^2 < 3p^k$$
contradicts
$$\frac{\sigma(m^2)}{p^k} \geq 7,$$
as the latter implies that
$$\frac{7p^k}{2} < m^2.$$

This completes the proof of Theorem \ref{VeryLargeLowerBound}.

\subsection{A proof of Corollary \ref{VeryLargeLowerBound2}}\label{Proof2}

Let $N = p^k m^2$ be an odd perfect number with special prime $p$.  By Lemma \ref{BroughanDelbourgoZhou}, we have
$$i(p)=\sigma(m^2)/p^k \geq 315.$$
This implies that
$$p^k < \frac{2m^2}{315}.$$
Adding $313{m^2}/315$ to both sides of the last inequality, and subtracting $p^k$, we get
$$m^2 - p^k > \frac{313m^2}{315}.$$
This finishes the proof of Corollary \ref{VeryLargeLowerBound2}.

\begin{remark}\label{NumericalLowerBound}
Note that we obtain, per Lemma \ref{estimateform}, the numerical lower bound
$$m^2 - p^k > 2m > 2\cdot{10}^{375}$$
from Theorem \ref{VeryLargeLowerBound}, and the numerical lower bound
$$m^2 - p^k > \frac{313m^2}{315} > \frac{313}{315}\cdot{10}^{750}$$
from Corollary \ref{VeryLargeLowerBound2}.
\end{remark}

\section{On Theorem \ref{MainResult1}}\label{Proof3}

\subsection{Sample Proof Arguments for $m<p^k$}
First, we look at the following sample proof arguments, by considering small values for $m^2 - p^k$.  (Lemma \ref{SampleProofArgument1} is exactly the same as Theorem \ref{PreludeToMainResult} in Section \ref{Intro}, and is reproved here for context.)

\begin{lemma}\label{SampleProofArgument1}
If $N_1 = p^k m^2$ is an odd perfect number given in Eulerian form and satisfying $m^2 - p^k = 8$, then the inequality $m < p^k$ holds.
\end{lemma}

\begin{proof}
Let $N_1 = p^k m^2$ be an odd perfect number with special prime $p$, satisfying $$m^2 - p^k = 8.$$  Subtracting $9$ from both sides and transferring $p^k$ to the other side of the equation, we obtain
$$(m + 3)(m - 3) = m^2 - 9 = p^k - 1.$$
By Lemma \ref{estimateform}, we have $m > {10}^{375}$.  Also, trivially we know that $p^k \geq 5$.  Hence both LHS and RHS of the last equation are positive.

Since $m - 3$ is a positive integer, this implies that
$$(m + 3) \mid (p^k - 1)$$
from which we obtain
$$m < m + 3 \leq p^k - 1 < p^k.$$
\end{proof}

\begin{lemma}\label{SampleProofArgument2}
If $N_2 = p^k m^2$ is an odd perfect number given in Eulerian form and satisfying $m^2 - p^k = 40$, then the inequality $m < p^k$ holds.
\end{lemma}

\begin{proof}
Let $N_2 = p^k m^2$ be an odd perfect number with special prime $p$, satisfying $$m^2 - p^k = 40.$$  Subtracting $49$ from both sides and transferring $p^k$ to the other side of the equation, we obtain
$$(m + 7)(m - 7) = m^2 - 49 = p^k - 9.$$
By Lemma \ref{estimateform}, we have $m > {10}^{375}$.  Hence, the LHS, and therefore the RHS, of the last equation are positive.

Since $m - 7$ is a positive integer, this implies that
$$(m + 7) \mid (p^k - 9)$$
from which we obtain
$$m < m + 7 \leq p^k - 9 < p^k.$$
\end{proof}

Note that in the proofs of both Lemma \ref{SampleProofArgument1} and Lemma \ref{SampleProofArgument2}, we subtracted the \emph{nearest square that is larger than} the value of $m^2 - p^k$.

\subsection{A proof of Theorem \ref{MainResult1}} 
Let us now prove Theorem \ref{MainResult1}.

\begin{proof}
Lemma \ref{SampleProofArgument1} and Lemma \ref{SampleProofArgument2} can be easily generalized, as follows: Any positive integer $m^2-p^k=2^r t=s$ where $s$ is less than $(m-1)^2=m^2-2m+1$ will lead to $m<p^k$. Because we are sure that $(m-1)^2$ is the biggest possible square less than $m^2$, then any positive integer $s$ falling below $(m-1)^2$ will lead to $m<p^k$. So there is nothing special about $8$ or $40$, and other equally small integers will still work because they are all less than $(m-1)^2$. In particular, note by Lemma \ref{estimateform} that $m>{10}^{375}$ which means that 
$$(m-1)^2=m^2-2m+1=m(m-2)+1>({10}^{375})\cdot({{10}^{375}}-2)+1.$$
Since $m^2$ is so huge, therefore $(m-1)^2=m^2-2m+1$ is also extremely big and there are plenty of numbers that fall below it. Essentially the equation for any positive integer $s=2^r t$ which is less than $(m-1)^2=m^2-2m+1$ becomes:
$$m^2-(m-1)^2=p^k+s-(m-1)^2$$
Since $s<(m-1)^2$, therefore $s-(m-1)^2$ will be a negative integer [let us call it $-u$ (where $u>0$)] and the above equation can be rewritten as follows:
$$m^2-(m-1)^2=p^k-u$$
$$\bigl(m+(m-1)\bigr)\cdot\bigl(m-(m-1)\bigr)=p^k-u$$
This equation, which has been generalized, resembles the equations in Lemma \ref{SampleProofArgument1} and Lemma \ref{SampleProofArgument2}.
Final inequality becomes:
$$m<(m+(m-1))=p^k-u<p^k$$ 
This does not in any way automatically mean that $m<p^k$, since we still do not know whether $p^k$ is indeed greater than $2m-1$ or not. This is still an open question.
\end{proof}

We shall use a similar, yet slightly different technique in proving Theorem \ref{MainResult2}, subject to some minimal conditions that we will impose.

\section{On Theorem \ref{MainResult2}}\label{Proof4}

By Theorem \ref{QuantityIsNotASquare}, the quantity $m^2 - p^k$ is not a square.  Additionally, we know that $m^2 - p^k \equiv 0 \pmod 4$.

Thus, in general, we may write
$$m^2 - p^k = 2^r t$$
where we know that $2^r \neq t$, $r \geq 2$, and $\gcd(2,t)=1$. (Now, since $r \geq 2$ holds, it follows that $2^r t$ is not squarefree.  Theorem \ref{QuantityIsNotASquare} then implies that $2^r t$ can be expressed as the product of a square and a squarefree integer.)

Note that it is easy to prove the following lemmas.

\begin{lemma}\label{mnotequaltotwotother}
If $N = p^k m^2$ is an odd perfect number given in Eulerian form and satisfying $m^2 - p^k = 2^r t$ where $r \geq 2$ and $\gcd(2,t)=1$, then $m \neq 2^r$.
\end{lemma}

\begin{proof}
The proof follows from the fact that $r \geq 2$ and $m$ is odd.
\end{proof}

\begin{lemma}\label{mnotequaltot}
If $N = p^k m^2$ is an odd perfect number given in Eulerian form and satisfying $m^2 - p^k = 2^r t$ where $r \geq 2$ and $\gcd(2,t)=1$, then $m \neq t$.
\end{lemma}

\begin{proof}
Assume to the contrary that $p^k m^2$ is an odd perfect number with special prime $p$ satisfying $m^2 - p^k = 2^r t$ and $m = t$.

We obtain
$$m^2 - p^k = 2^r m$$
$$m^2 - 2^r m = p^k$$
$$m(m - 2^r) = p^k.$$
Since $p^k \geq 5$, it follows that $m - 2^r > 0$, and therefore that $m \mid p^k$.  

This last divisibility constraint contradicts $\gcd(p,m)=1$.
\end{proof}

\subsection{A proof of Theorem \ref{MainResult2}}
We are now ready to prove Theorem \ref{MainResult2}.

\begin{proof}
Six cases need to be considered from Lemma \ref{mnotequaltotwotother}, Lemma \ref{mnotequaltot}, and the constraint $2^r \neq t$ (which holds by Theorem \ref{QuantityIsNotASquare}):
\begin{enumerate}
\item{$m > t > 2^r$}
\item{$m > 2^r > t$}
\item{$t > m > 2^r$}
\item{$2^r > m > t$}
\item{$t > 2^r > m$}
\item{$2^r > t > m$}
\end{enumerate}

We consider these six cases in turn below:

Case (1): $m > t > 2^r$

Note that Case (1) implies $m - t > 0$ and 
$$2^{2r} < m^2 - p^k = 2^r t < t^2.$$
Following our method, we subtract $t^2$ from both sides of $m^2 - p^k = 2^r t$ to obtain
$$(m + t)(m - t) = m^2 - t^2 = p^k - t(t - 2^r).$$
Since $m - t$ is a positive integer, both sides of the last equation are positive. This then implies that
$$(m + t) \mid \bigg(p^k - t(t - 2^r)\bigg)$$
from which we obtain
$$m < m + t \leq p^k - t(t - 2^r) < p^k,$$
since $t > 2^r$.

Case (2): $m > 2^r > t$

Note that Case (2) implies $m - 2^r > 0$ and
$$t^2 < m^2 - p^k = 2^r t < 2^{2r}.$$
Following our method, we subtract $2^{2r}$ from both sides of $m^2 - p^k = 2^r t$ to obtain
$$(m + 2^r)(m - 2^r) = m^2 - 2^{2r} = p^k - 2^r (2^r - t).$$
Since $m - 2^r$ is a positive integer, both sides of the last equation are positive. This then implies that
$$(m + 2^r) \mid \bigg(p^k - 2^r (2^r - t)\bigg)$$
from which we obtain
$$m < m + 2^r \leq p^k - 2^r (2^r - t) < p^k,$$
since $2^r > t$.

Case (3): $t > m > 2^r$

Note that Case (3) implies $(m - t)(m - 2^r) < 0$, which implies that
$$m^2 + 2^r t < m(2^r + t)$$
$$p^k = m^2 - 2^r t < m^2 + 2^r t < m(2^r + t),$$
from which we cannot conclude whether $p^k < m$ or $m < p^k$.

On the other hand, the inequality
$$m^2 + 2^r t < m(2^r + t)$$
may be rewritten as
$$m^2 + (m^2 - p^k) < m(2^r + t)$$
$$2m^2 < m(2^r + t) + p^k.$$
Since we want to prove $m < p^k$, assume to the contrary that $p^k < m$.  We obtain
$$2m^2 < m(2^r + t) + p^k < m(2^r + t) + m = m(2^r + t + 1),$$
from which it follows that
$$2m < 2^r + t + 1.$$

It may be possible to derive a contradiction from this last inequality under this case, by considering the estimate in Lemma \ref{estimateform}.

The authors leave Case (3) as an open problem for other researchers to investigate.

Case (4): $2^r > m > t$

Note that Case (4) implies $(m - t)(m - 2^r) < 0$, which implies that
$$m^2 + 2^r t < m(2^r + t)$$
$$p^k = m^2 - 2^r t < m^2 + 2^r t < m(2^r + t),$$
from which we cannot conclude whether $p^k < m$ or $m < p^k$.

On the other hand, the inequality
$$m^2 + 2^r t < m(2^r + t)$$
may be rewritten as
$$m^2 + (m^2 - p^k) < m(2^r + t)$$
$$2m^2 < m(2^r + t) + p^k.$$
Since we want to prove $m < p^k$, assume to the contrary that $p^k < m$.  We obtain
$$2m^2 < m(2^r + t) + p^k < m(2^r + t) + m = m(2^r + t + 1),$$
from which it follows that
$$2m < 2^r + t + 1.$$

It may be possible to derive a contradiction from this last inequality under this case, by considering the estimate in Lemma \ref{estimateform}.

The authors leave Case (4) as an open problem for other researchers to investigate.

Case (5): $t > 2^r > m$

Note that Case (5) implies $m < t$ and $m < 2^r$, which means that $m^2 < 2^r t$.  Thus, $$m^2 - 2^r t < 0.$$
This contradicts $m^2 - 2^r t = p^k \geq 5$.  Hence, Case (5) does not hold.

Case (6): $2^r > t > m$

Note that Case (6) implies $m < t$ and $m < 2^r$, which means that $m^2 < 2^r t$.  Thus, $$m^2 - 2^r t < 0.$$
This contradicts $m^2 - 2^r t = p^k \geq 5$.  Hence, Case (6) does not hold.

This concludes the proof of Theorem \ref{MainResult2}.
\end{proof}

\section{Concluding remarks}\label{Conclusion}

The first-named author, together with Dagal, initially attempted an unconditional proof for $m < p^k$ in November 2020 \cite{DagalDris}.  Several errors, however, were identified by Ochem and the anonymous user mathlove in MathOverflow [\url{https://mathoverflow.net/questions/376268}].  Ochem pointed out that the condition $$0 < m-\lceil\sqrt{m^2-p^k}\rceil$$ requires $p^k \geq 2m-1$, which would be an unhelpful assumption since the goal is to prove $p^k > m$.  (Note that we are sure that $p^k \neq 2m - 1$, because otherwise the quantity $m^2 - p^k = m^2 - 2m + 1 = (m - 1)^2$ would be a square, contradicting Theorem \ref{QuantityIsNotASquare}.)  The authors emphasize that, while the method used is similar, this paper does not substantially rely on nor use any crucial results from the preprint \cite{DagalDris}.

This paper is an attempt at resolving the difficulties in that earlier preprint \cite{DagalDris}, carefully delineating the particular cases that need to be considered (i.e. Cases (1) through (6) in Section \ref{Proof3}).

\subsection{Future Research}\label{FutureResearch}

Indeed, the following cases remain to be considered, which the authors leave as open problems for other researchers to investigate:
\begin{enumerate}
\item{$t > m > 2^r$}
\item{$2^r > m > t$}
\end{enumerate}
Here, $N = p^k m^2$ is an odd perfect number with special prime $p$ satisfying $m^2 - p^k = 2^r t$ where $r \geq 2$ and $\gcd(2,t)=1$.

It might be prudent to note the following simpler proof (taken from [\url{https://mathoverflow.net/questions/406701}]) for Theorem \ref{MainResult2}: We consider the same Cases (1) through (6).  Case (5): $t > 2^r > m$ and Case (6): $2^r > t > m$ are easily dispensed with.  Under Case (1) and Case (2), we can prove that the inequality $m < p^k$ holds, as follows:

Under Case (1): $m > t > 2^r$, we have:
$$(m - t)(m + 2^r) > 0$$
$$p^k = m^2 - 2^r t > m(t - 2^r) = m\left|2^r - t\right| \geq m.$$

Under Case (2): $m > 2^r > t$, we have:
$$(m - 2^r)(m + t) > 0$$
$$p^k = m^2 - 2^r t > m(2^r - t) = m\left|2^r - t\right| \geq m.$$

So we are now left with Case (3) and Case (4).

Under Case (3): $t > m > 2^r$, we have:
$$(m + 2^r)(m - t) < 0$$
$$p^k = m^2 - 2^r t < m(t - 2^r) = m\left|2^r - t\right|.$$

Under Case (4): $2^r > m > t$, we have:
$$(m - 2^r)(m + t) < 0$$
$$p^k = m^2 - 2^r t < m(2^r - t) = m\left|2^r - t\right|.$$

Note that, under Case (3) and Case (4), we actually have
$$\min(2^r,t) < m < \max(2^r,t).$$

But the condition $\left|2^r - t\right|=1$ is sufficient for $p^k < m$ to hold, under Case (3) and Case (4).  (Note that the condition $\left|2^r - t\right|=1$ is not necessary for $p^k < m$ to hold, under Case (3) and Case (4), basically because $\min(2^r,t) < m < \max(2^r,t)$ means that $2^r$ and $t$ are not consecutive integers.)

Therefore, if one could show that the condition $\left|2^r - t\right|=1$ holds unconditionally, then one would be able to prove that $m < p^k$ (since we would then know that Case (3) and Case (4) could not occur).

It is therefore natural to predict the truth of the following conjecture:

\begin{conjecture}\label{ConjectureOPNLinkToEPN}
If $N = p^k m^2$ is an odd perfect number given in Eulerian form and satisfying $m^2 - p^k = 2^r t$, then $t = 2^r - 1$ is prime.
\end{conjecture}

\subsection{Final Remarks}
Note that we may summarize Case (3) and Case (4) as $\min(2^r,t) < m < \max(2^r,t)$, and Case (1) and Case (2) as $m > \max(2^r, t)$.  Lastly, note that Case (5) and Case (6) can be summarized as $m < \min(2^r,t)$.  Notice further that we have
$$\min(2^r,t)\cdot\max(2^r,t) = 2^r t$$
and that
$$\max(2^r,t)-\min(2^r,t)=\bigg(\frac{(2^r+t)+\left|2^r - t\right|}{2}\bigg)-\bigg(\frac{(2^r+t)-\left|2^r - t\right|}{2}\bigg)$$
$$=\frac{2\cdot\left|2^r - t\right|}{2}=\left|2^r - t\right|.$$

First, consider Case (5) and Case (6).  We obtain the series of implications
$$m < \min(2^r, t) \Rightarrow \frac{1}{\min(2^r, t)} < \frac{1}{m} \Rightarrow \max(2^r, t) = \frac{2^r t}{\max(2^r, t)} < \frac{m^2 - p^k}{m} = m-\frac{p^k}{m} < m,$$
resulting in the contradiction
$$\max(2^r, t) < m < \min(2^r, t).$$
Hence, Case (5) and Case (6) do not hold, which confirms what was previously proved in this regard.

Consider Case (1) and Case (2).  We obtain the series of implications
$$m > \max(2^r, t) \Rightarrow \frac{1}{\max(2^r, t)} > \frac{1}{m} \Rightarrow \min(2^r, t) = \frac{2^r t}{\max(2^r, t)} > \frac{m^2 - p^k}{m}.$$
But from Corollary \ref{VeryLargeLowerBound2}, we have the lower bound
$$m^2 - p^k > \frac{313m^2}{315}.$$
We then have
$$\min(2^r, t) > \frac{m^2 - p^k}{m} > \frac{313m}{315}$$
so that we have
$$\max(2^r, t) < m < \frac{315}{313}\cdot\min(2^r,t).$$
(No contradictions, thus far.  Note that $\max(2^r, t) < \min(2^r, t)$ would have been a contradiction.)

Now, consider Case (3) and Case (4).  We obtain the series of implications
$$\min(2^r, t) < m < \max(2^r, t) \Rightarrow \frac{1}{\max(2^r, t)} < \frac{1}{m} < \frac{1}{\min(2^r, t)}$$ 
$$\Rightarrow \min(2^r, t) = \frac{2^r t}{\max(2^r, t)} < \frac{m^2 - p^k}{m} < \frac{2^r t}{\min(2^r, t)} = \max(2^r, t).$$
(No contradictions thus far here, as well.)

These final remarks illustrate some of the difficulties in trying to arrive at a contradiction for either Case (1), Case (2), Case (3), or Case (4).

\section*{Acknowledgments} I.~T.~S.~D. thanks TUA-URDC for support.  We thank the anonymous referees who have made invaluable suggestions for improving the quality of this paper.  J.~A.~B.~D. likewise thanks the anonymous MathOverflow user mathlove [\url{https://mathoverflow.net/users/34490}] and Tony Kuria Kimani for sharing their expertise.

\end{document}